\documentclass[12pt,reqno]{amsart}

\usepackage{graphicx}

\usepackage{amssymb}
\usepackage[inline]{enumitem}
\usepackage{tikz,lipsum,lmodern}
\usepackage[most]{tcolorbox}
 \usepackage{youngtab}
 \usepackage{ytableau}
 \usepackage{hyperref}
\usepackage{amsthm}
\theoremstyle{plain}
\usepackage{xcolor}
\usepackage{tikz}
\usetikzlibrary{
  graphs,
  graphs.standard,
  arrows.meta
}

\definecolor{myblue}{RGB}{80,80,160}
\definecolor{mygreen}{RGB}{80,160,80}

\newtheorem*{theorem*}{Theorem}
\usepackage{lineno}

\title{Tiling Ferrers Diagrams Handout}

\author{David Jiang}

\begin{document}
\begin{abstract}
We will show that a necessary and sufficient condition for a Ferrers board to be fully tileable with $1 \times 2$ dominoes requires the board to be $2-$colorable such that no color is adjacent to its own color using both induction and a graph theory approach. We will walk through all prerequisite knowledge and go through the failed attempts we tried while also providing supplementary exercises that fit the topics. Plenty of background content is included, so even if you don't know much about the subject, it should still be readable. If you do know most of the background content, feel free to skip around.
\end{abstract}

\maketitle
\tableofcontents
\newpage
\section{\textbf{Partitions}} In order to get started, we're first going to lay the ground work for the first few mathematical objects that we will be working with. Starting with just partitions.
\begin{tcolorbox}
[colback=blue!5!white,colframe=blue!75!black,title=\textbf{Definition 1.1} (Integer Partitions)]
An integer partition is a way of writing a positive integer as a sum of other positive integers, where the order of the summands does not matter. 
\end{tcolorbox}

\begin{tcolorbox}[colback=yellow!5!white,colframe=yellow!75!black,title=\textbf{Example 1.1} Partitions of 5]
  \begin{itemize}
        \item $5 = 1+1+1+1+1$
        \item $5 = 2+1+1+1$
        \item $5 = 2+2+1$
        \item $5 = 3+1+1$
        \item $5 = 3+2$
        \item $5 = 4+1$
        \item $5 = 5$
    \end{itemize}
\end{tcolorbox}
\begin{tcolorbox}[colback=blue!5!white,colframe=blue!75!black,title=\textbf{Definition 1.2} (Ferrers Boards) ]
  The Ferrers Board is a pictorial representation of an integer partition. Taking a certain partition, the Ferrers Board is generated by taking the largest part $n_1$ of an integer partition and placing $n_1$ many blocks on the first row, for the second largest part $n_2$ placing $n_2$ many blocks on the second row, and continuing until all parts of the partition have been represented with blocks.
\end{tcolorbox}

\begin{tcolorbox}[colback=yellow!5!white,colframe=yellow!75!black,title=\textbf{Example 1.2} Ferrers Boards of 5]
  $\young(~,~,~,~,~)$ \hspace{55pt}
    $\young(~~,~,~,~)$ \hspace{55pt}
    $\young(~~,~~,~)$ \hspace{55pt}
    $\young(~~~,~,~)$ \hspace{55pt}\\\\
    $\young(~~~,~~)$ \hspace{28pt}
    $\young(~~~~,~)$ \hspace{25pt}
    $\young(~~~~~~)$
\end{tcolorbox}

\begin{tcolorbox}[enhanced,attach boxed title to top center={yshift=-3mm,yshifttext=-1mm},
  colback=blue!5!white,colframe=blue!75!black,colbacktitle=red!80!black,
  title=Exercise 1.x,fonttitle=\bfseries,
  boxed title style={size=small,colframe=red!50!black} ]
  \textbf{1.1 ($\star$)} List and draw the Ferrers board for the partitions of 6.
    \\\\
    \textbf{1.2 ($\star \star$)} Show that the partitions can be upper bounded by the Fibonacci numbers.\\\\
    \textbf{1.3 ($\star \star$)} Show that partitions into distinct odd parts is equal to partitions that have symmetric Ferrers diagrams when you reflect it. (Hint: look at the Ferrer diagrams, what parts must be odd?)
\end{tcolorbox}
\newpage
\section{\textbf{Domino Tilings}} Building off of the classic kids toys of dominoes, they are $1\times 2$ blocks. We will be looking at grids and see if we can find heuristic methods to prove that certain grids are tilable and find methods to generate such tilings.\\\\
We will first look at some examples to get an idea of what to do. For the sake of simplicity, all dominoes refers to $1 \times 2$ dominoes. More complicated dominoes do exist but we will not be discussing them here.

\begin{tcolorbox}[enhanced,attach boxed title to top center={yshift=-3mm,yshifttext=-1mm},
  colback=blue!5!white,colframe=blue!75!black,colbacktitle=red!80!black,
  title=Exercise 2.x,fonttitle=\bfseries,
  boxed title style={size=small,colframe=red!50!black} ]
  \textbf{2.1 ($\star$)} 
  Find all domino tilings of \hspace{5pt} $\young(~~~,~~~)$
\\\\
  \textbf{2.2 ($\star$)} Show that it is possible to tile all $2n \times m$ boards where $n,m$ are positive integers. 
\\\\
  \textbf{2.3 ($\star\star$)} Given an $8 \times 8$ board, what 2 corners can you remove and still be able to tile the board?
\\\\
    \textbf{2.4 ($\star\star$)} Show that a stair case grid can't be tiled with dominoes.\\
    \textbf{Example:} $$\young(~~~~,~~~,~~,~)$$
\\\\
\textbf{2.5 ($\star\star\star$)} Given a $2n \times m$ checkboard where $n,m$ are integers. Show that when you remove a white and a black square from the board, the board will still remain tileable.\footnote[1]{This is called: Gomory's theorem}
\end{tcolorbox}

Now that we have an idea of what domino tilings and Ferrers Diagrams are, we can introduce the question at hand. Given a Ferrers diagram, how can we look at it and determine if it's able to be tiled or not? Or more specifically:
\begin{tcolorbox}
[colback=blue!5!white,colframe=blue!75!black,title=\textbf{Question.}]
Find both a necessary and sufficient condition for when a Ferrers Board can be tiled by dominoes.
\end{tcolorbox}
\begin{tcolorbox}[enhanced,attach boxed title to top center={yshift=-3mm,yshifttext=-1mm},
  colback=blue!5!white,colframe=blue!75!black,colbacktitle=red!80!black,
  title=Exercise 3.x,fonttitle=\bfseries,
  boxed title style={size=small,colframe=red!50!black} ]
\textbf{2.6 ($\star$)} Find a tiling for the following Ferrers boards:
\begin{enumerate}
    \item $\young(~~~~,~~,~,~)$
    \item $\young(~~~,~~~,~~~,~)$
    \item $\young(~~~~~,~~~~~,~~~,~)$
\end{enumerate}
\end{tcolorbox}
We showed in exercise 2.4 that a staircase grid can't be tiled. Since the grid is essentially a Ferrers diagram, we can make the following lemma.
\begin{tcolorbox}
[colback=green!5!white,colframe=green!75!black,title=\textbf{Lemma 1.}]
A Ferrers board in which each sequential part differs by only $1$ can't be tiled.\\\\
\textbf{Example}
$$\young(~~~~~~,~~~~~,~~~~)$$
$$\vdots$$
Can't be tiled
\end{tcolorbox}
\newpage
\section{\textbf{Inductive Algorithm}}
This section is a bit drawn out because it details all the observations we made and attempts we tried. Most is useless information, skip to the algorithm if you want.
\subsection{Necessary and Sufficient Conditions}
Before we describe the algorithm, it is important to go over what ideas we have for necessary and sufficient conditions. Depending on how you solved 2.5, you could use a \href{https://en.wikipedia.org/wiki/Hamiltonian_path}{Hamiltonian Path} argument. An even length  Hamiltonian Path is indeed a sufficient condition, because if you have such path, you can lay dominoes along the path and a tiling will be formed.
\begin{tcolorbox}[enhanced,attach boxed title to top center={yshift=-3mm,yshifttext=-1mm},
  colback=blue!5!white,colframe=blue!75!black,colbacktitle=red!80!black,
  title=Exercise 3.x,fonttitle=\bfseries,
  boxed title style={size=small,colframe=red!50!black} ]
\textbf{3.1 ($\star\star$)} Find a Ferrers board that can be tiled, but does not contain a Hamiltonian path.
\end{tcolorbox}
Implied by the exercise, a Hamiltonian path is not a necessary and sufficient condition. Let's think about another previous exercise. Depending on how you solved 2.3, you could consider coloring the Ferrers diagram. One such solution for that question argues that if you color the grid with black and white alternating squares like a checkerboard, then in order for a tiling to be possible, you must have an even amount of equal black and white squares. This is because each domino can only cover one black and one white square, so if you have an unequal amount of each, a tiling is not possible no matter how hard you try. We can use this argument for our Ferrers boards. Let's imagine superimposing our Ferrers board onto a checkerboard. We know that if a Ferrers board can be tiled, it would imply that there is an equal amount of even black and white squares on the board. So it is a necessary condition, is it a sufficient condition?

\begin{tcolorbox}[enhanced,attach boxed title to top center={yshift=-3mm,yshifttext=-1mm},
  colback=blue!5!white,colframe=blue!75!black,colbacktitle=red!80!black,
  title=Exercise 3.x,fonttitle=\bfseries,
  boxed title style={size=small,colframe=red!50!black} ]
\textbf{3.3 ($\star\star\star\star\star$)} Find a Ferrers board that can't be tiled but contains an equal amount of even white and black tiles. (Hint: Keep reading)
\end{tcolorbox}
That may have been evil if you tried that exercise because this is the necessary and sufficient condition we use and prove.

\subsection{Algorithm}
A common approach to problems such as this question is finding an algorithm that completely tiles a Ferrers board. This is because if you can prove an algorithm to always tile a board then that would imply necessary and sufficient conditions for it to be tiled.\\\\
Proving an algorithm requires two things primarily. Proving that the algorithm terminates and proving that the base cases are satisfied.\\\\
Termination is hard to prove without an actual algorithm in mind, so let's first think about what potential base cases we could have.\\\\
We showed in exercise 2.2 that if we can reduce something to a $2n \times m$ grid then it can be tileable \footnote[1]{If you are interested in these kinds of ideas, this is based on the idea of Durfee squares. A lot of identites come from Durfee squares.}. Perhaps we can tile the largest rectangle contained within a Ferrers board first and then inductively work on the rest of the non-rectangle parts.\\
\begin{tcolorbox}[enhanced,attach boxed title to top center={yshift=-3mm,yshifttext=-1mm},
  colback=blue!5!white,colframe=blue!75!black,colbacktitle=red!80!black,
  title=Exercise 3.x,fonttitle=\bfseries,
  boxed title style={size=small,colframe=red!50!black} ]
\textbf{3.1 ($\star$)} Find an example where when you fill in the largest $2n\times m$ rectangle first in a Ferrer diagram, the left over parts won't be tileable.
\end{tcolorbox} So this heuristic doesn't quite work. And more importantly, it doesn't use our necessary and sufficient condition in any way. Let's look back at another exercise, depending on how you solved 1.3, you might notice that L shapes with even amount of tiles in the Ferrer board always can be tiled. For example:
$$\young(\bullet\bullet\bullet,\bullet~~,\bullet~~,\bullet)$$
So we can tile this L shape and then reduce it down to a smaller Ferrers diagram that maintains the white/black property of a supposedly tileable board. Again, there are issues with this method. What happens if you have an odd shaped L? Now we will detail the algorithm we designed that tiles every board with the white/black property, as well as prove it afterwards.
\subsection{Actual Algorithm (From Jonah Guse)}
\subsubsection{\textbf{General Idea:}} We imagine superimposing a Ferrers board onto a chess board, so that a square with coordinates $(i, j)$ is colored black if $i + j$ is even, and colored white instead if $i + j$ is odd. When tiling such a colored Ferrers board with dominoes, each domino must cover exactly one black square and one white square. Therefore, if a Ferrers board can be tiled by dominoes, the number of black squares equals the number of white squares.\\\\
In fact, the converse is also true! If a Ferrers board has an equal number of black and white squares, then it can be tiled by dominoes. To prove this, we consider a recursive algorithm to cover such a Ferrers board with dominoes.
\subsubsection{\textbf{Algorithm}} \hspace{1pt}
\\ \textbf{Base Case:}
If our Ferrer board is empty, then we are done.\\\\
\textbf{Case 1 (Left):}
If the leftmost column contains an even number of squares, we can cover it with vertical dominoes, leaving a smaller (or empty) Ferrer board with an equal number of white and black squares. We tile this smaller board using this algorithm.\\\\
\textbf{Case 2 (Top):}
If the topmost row contains an even number of squares, we can cover it with horizontal dominoes, leaving a smaller (or empty) Ferrer board with an equal number of white and black squares. We tile this smaller board using this algorithm.\\\\
\textbf{Case 3 (Bottom):}
If the lowest nonempty row contains more than 1 square, we cover the rightmost two squares in this row with a domino, leaving a smaller Ferrer board with an equal number of white and black squares. We tile this smaller board using this algorithm.\\\\
\textbf{Case 4 (Odd):}
If all 3 cases are not true then we have case 4. Without loss of generality let the top left square be black. Then the number of black squares in the first row is one greater than the number of white squares in the first row. Furthermore, the only square in the bottom row is black.
\\\\
We now label the kth row with the number of black squares in the first k rows minus the number of white squares in the first k rows. Thus, the first row is labeled with the number 1. Furthermore, the last row is labeled with the number 0 and contains exactly one black square, so the second to last row is labeled with the number $-1$.\\\\
Because the squares are colored in an alternating pattern in each row, the label on any row after the first differs from the label above it by exactly 0 or 1. Therefore, because the first row is labeled with the number 1, and the second to last row is labeled with the number -1, some variation of intermediate value theorem tells us that some row in between is labeled with the number 0.
\\\\
Then we choose such a row labeled with a 0. This row and all rows above it then form a smaller Ferrer board with an equal number of white and black squares. Cover this smaller board using this algorithm. Similarly, all rows below this row also form a smaller Ferrer board with an equal number of white and black squares. Cover this smaller board using this algorithm.\\\\

\begin{tcolorbox}[colback=yellow!5!white,colframe=yellow!75!black,title=\textbf{Example 3} Algorithm on a Tiling]
  $$\young(~~~~~~~~,~~~~~~,~~~~~,~~~~,~~~~,~)$$
  Since the left most column has $6$ tiles that means we have case 1, we can first tile that up.
    
$$\young(\bullet~~~~~~~,\circ~~~~~,\bullet~~~~,\circ~~~,\bullet~~~,\circ)$$ 
This gives us a new Ferrer board that maintains the white/blackness property that we are looking for.\\\\
Now we see that both the top column and left row are odd so we look at our next condition. We notice that the bottom row has more than $1$ tile which is case 3, so we can tile that as well, which gives us.
$$\young(\bullet~~~~~~~,\circ~~~~~,\bullet~~~~,\circ~~~,\bullet~\circ\bullet,\circ)$$
Now again we have a smaller Ferrers board that maintains our target property. Since we didn't change the left row or top column, this means we need to look at case 4.
\begin{center}
\includegraphics[width=50mm]{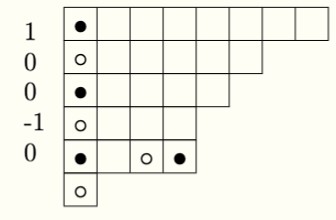}
\end{center}
So by this label, this means we can split it into smaller Ferrers Diagrams. Particularly we can split it into:\\
\end{tcolorbox}
\begin{tcolorbox}[colback=yellow!5!white,colframe=yellow!75!black,title=\textbf{Example 3} Algorithm on a Tiling]
$\young(~~~~~~~,~~~~~)$ \hspace{55pt}
    $\young(~~~~)$ \hspace{55pt}
    $\young(~~~,~)$ \hspace{55pt}\\
In each one of these we can apply the algorithm again. Let's first look at the first one. The left most column now has an even amount of tiles so we can tile it. We can actually do this for most of the columns because we have a $5\times 2$ rectangle contained in the diagram.
\begin{center}
$\young(\bullet\circ\bullet\circ\bullet~~,\circ\bullet\circ\bullet\circ)$
\end{center}
Then we notice that the remaining part can be tiled with a singular tile so this part of the board is tileable. The other 2 parts are also easily tileable. So this gives us the final tiling:
$$\young(\bullet\circ\bullet\circ\bullet\circ\bullet\bullet,\bullet\circ\bullet\circ\bullet\circ,\circ\bullet\bullet\circ\circ,\circ\circ\bullet\bullet,\bullet\circ\circ\circ,\bullet)$$
\end{tcolorbox}

\subsection{\textbf{Proof}} As mentioned earlier, the two parts of an algorithmic proof is to show that it terminates and has a base case. In this case the base case is already outlined in the algorithm so it exists. Then we also need to prove termination of the algorithm.

\begin{tcolorbox}[enhanced,attach boxed title to top center={yshift=-3mm,yshifttext=-1mm}, colback=blue!5!white,colframe=blue!75!black,colbacktitle=red!80!black,
  title=Exercise 3.x,fonttitle=\bfseries,
  boxed title style={size=small,colframe=red!50!black} ]
\textbf{3.2 ($\star\star$)} Show that the algorithm described above terminates.
\end{tcolorbox}
\newpage
\section{\textbf{Graph Theory Solution}}
We've shown that our solution can be solved using an inductive algorithm. Now we will show a graph theoretic approach. Before we do so, we will detail some topics in graph theory to make sure that the solution is clear.

\subsection{Introduction To Bipartite Graphs}
We will primarily be working with bipartite graphs. 
\begin{tcolorbox}
[colback=blue!5!white,colframe=blue!75!black,title=\textbf{Definition 4.1} (Bipartite Graphs)]
Given a graph $G = (V,E)$, if $V$ can be split into 2 groups $V_1, V_2$ such that $V_1 \cap V_2 = \varnothing$ so no two graph vertices within the same set are adjacent. And where each element in $V_1$ is adjacent to some element in $V_2$. Then $G$ is a bipartite graph.
\end{tcolorbox}
\begin{tcolorbox}[colback=yellow!5!white,colframe=yellow!75!black,title=\textbf{Example 4.1} Identifying Bipartite Graphs]
\begin{center}
  \begin{tikzpicture}
   \graph[nodes={draw, circle,fill=myblue}, radius=.5cm,
           empty nodes, branch down=1 cm,
           grow right sep=4cm] {subgraph I_nm [V={a, b, c, d, e}, W={1,...,4}];
  a -- { 1};
  b -- { 1, 2 };
  c -- { 2,3, 4 };
  d -- {1,2};
  e -- { 1,4}
};
\end{tikzpicture}
\end{center}
The nodes on the left can be $V_1$ and the nodes on the right are $V_2$. We observe that $V_1 \cap V_2 = \varnothing$, and that each element is attached to another element in the other set.
\end{tcolorbox}
The reason we are talking about this, is because we can represent our Ferrers board as a bipartite graph. Since when we superimpose our board with a checkerboard, we can separate the squares into 2 separate disjoint subsets by partitioning them based on if its a white square or a black square. Since each square can only be connected to a square of opposite color, we see that this indeed is a bipartite graph.
\subsection{Perfect Bipartite Matching}
It might not be immediately obvious why we want to represent the Ferrers board as a bipartite graph. The main reason is because we can find a bipartite matching on a bipartite graph which represent a certain tiling of a board. If we find a perfect bipartite matching then it correlates with a perfect tiling of the Ferrer diagram.

\begin{tcolorbox}
[colback=blue!5!white,colframe=blue!75!black,title=\textbf{Definition 4.2} (Perfect Bipartite Matching)]
A perfect bipartite matching is where every vertex in the left partition of a bipartite graph is matched to a vertex in the right partition, and vice versa. In other words, a perfect bipartite matching is a set of edges in a bipartite graph such that every vertex in the graph is incident to exactly one edge in the matching.
\end{tcolorbox}

\begin{tcolorbox}[colback=yellow!5!white,colframe=yellow!75!black,title=\textbf{Example 4.2} Perfect Bipartite Matching]
\begin{center}
  \begin{tikzpicture}
   \graph[nodes={draw, circle,fill=myblue}, radius=.5cm,
           empty nodes, branch down=1 cm,
           grow right sep=4cm] {subgraph I_nm [V={a, b, c, d}, W={1,...,4}];
  a -- {1};
  a -- [{red,thick}] {3};
  b -- [{red,thick}]{1};
  b -- {2};
  c -- { 3, 4 };
  c -- [{red,thick}]{2};
  d -- {1};
  d -- [{red,thick}]{4}
};
\end{tikzpicture}
\end{center}
If you look at just the red edges, they form a perfect bipartite matching where all of the left nodes are connected to all of the right nodes. So now you have $4$ edges and $8$ nodes.
\begin{center}
  \begin{tikzpicture}
   \graph[nodes={draw, circle,fill=myblue}, radius=.5cm,
           empty nodes, branch down=1 cm,
           grow right sep=4cm] {subgraph I_nm [V={a, b, c, d}, W={1,...,4}];
  a -- [{red,thick}] {3};
  b -- [{red,thick}]{1};
  c -- [{red,thick}]{2};
  d -- [{red,thick}]{4}
};
\end{tikzpicture}
\end{center}
\end{tcolorbox}
\subsection{Hall's Theorem}
We will be citing Hall's theorem for our solution, in this section I will detail what Hall's theorem is, but because the proof is quite lengthy I will leave it as an exercise. this will provide a detailed solution for our problem. First I will try to motivate the theorem by yourself.
\begin{tcolorbox}[enhanced,attach boxed title to top center={yshift=-3mm,yshifttext=-1mm},
  colback=blue!5!white,colframe=blue!75!black,colbacktitle=red!80!black,
  title=Exercise 4.x,fonttitle=\bfseries,
  boxed title style={size=small,colframe=red!50!black} ]
\textbf{4.1 ($\star$)} Given the following bipartite graph shown previously, show that all of the following subgraphs have a perfect bipartite matching.
\begin{center}
  \begin{tikzpicture}
   \graph[nodes={draw, circle,fill=myblue}, radius=.5cm,
           empty nodes, branch down=1 cm,
           grow right sep=4cm] {subgraph I_nm [V={a, b, c, d}, W={1,...,4}];
  a -- {1,3};
  b -- {1,2};
  c -- { 2,3, 4 };
  d -- {1,4};
};
\end{tikzpicture}
\end{center}
\begin{center}
\rule{1\textwidth}{.4pt}
\end{center}
\begin{center}
\begin{tikzpicture}
   \graph[nodes={draw, circle,fill=myblue}, radius=.5cm,
           empty nodes, branch down=1 cm,
           grow right sep=4cm] {subgraph I_nm [V={a, b}, W={1,...,3}];
  a -- {1,3};
  b -- {1,2};
};
\end{tikzpicture}
\end{center}
\begin{center}
\rule{1\textwidth}{.4pt}
\end{center}
\begin{center}
  \begin{tikzpicture}
   \graph[nodes={draw, circle,fill=myblue}, radius=.5cm,
           empty nodes, branch down=1 cm,
           grow right sep=4cm] {subgraph I_nm [V={a, b, c, d}, W={1,...,4}];
  a -- {1,3};
  c -- { 2,3, 4 };
  d -- {1,4};
};
\end{tikzpicture}
\end{center}
\textbf{4.2 ($\star\star\star$)} Prove or disprove that any bipartite graph has a perfect bipartite matching if and only if all proper subgraphs of it have a perfect bipartite matching.
\end{tcolorbox}
The point of the exercises is to get you to think about Hall's theorem, and more specifically to our specific version that we need.
\begin{tcolorbox}
[colback=blue!5!white,colframe=blue!75!black,title=\textbf{Definition 4.3} (Hall's Theorem)]
A bipartite graph $G = (V, E)$, with the bipartition $V = L\cup R$ where
$|L| = |R| = n$, has a perfect matching if and only if for every subset $S \subset L$ , $|N(S)| \geq |S|$.
\\\\Where $N(S)$ denotes the neighborhood of $S$, which is the set of nodes that the subset $S$ gets mapped to in $G$.
\end{tcolorbox}
You might have an idea of where the solution is going, but before we get to that point, we're going to prove this theorem and talk about it specifically for our problem.\\
\begin{tcolorbox}
[colback=gray!5!white,colframe=gray!75!black,title=\textbf{Proof 4.3}] First we will prove the \textbf{forward} direction.\\\\
Given the graph $G = (V,E)$ that has the bipartition $V = L \cup R$ and $|L|=|R|=n$, let's assume towards a contradiction that there exists $S\subset L$ such that $|N(S)| < |S|$. This means that there are some vertices in $S$ that aren't connected to distinct vertices in $R$. However, then that means it won't be a perfect matching, and therefore all subsets must satisfy $|N(S)| \geq |S|$.\\\\ Now let's prove the \textbf{backwards}\footnote{Adapted from Duke CS} direction.\\\\
Let's construct a perfect matching with induction on $S$, given the assumption that $|N(S)| \geq |S|$ for all subsets $S \subset L$.\\\\
\textbf{Base Case:} For $|L| = 1$. Let's let $l$ be the only vertex in $L$. Since we are provided our assumption from the theorem, we know that $|N(L)| \geq 1$. That means that $l$ must have 1 or more nodes in $R$ that it connects to. Let $M = ((l,r))$ where $r$ is a neighbor of $l$, then this gives us a perfect matching of size $1$ because you can just connect $l$ to $r$. So every vertex in $L$ is matched uniquely to one in $R$\\\\
\textbf{Inductive Solution:} Let's first assume that for all values $1 \leq k \leq |L|-1$ has a perfect bipartite matching of size $|S| = k$ if it follows our theorem.\\\\
Now try to complete the rest of the proof.
\end{tcolorbox}

\subsection{Graph Theoretic Solution}
Creating a bipartite graph from a Ferrers board is not too difficult. Since we are already trying to prove that the black/whiteness of the board is a necessary and sufficient condition, let's work from there. We know that each white square can only be connected to adjacent black squares and vice versa. This means that we already have $2$ sets of disjoint vertices that connect to only each other. We can denote all the white vertices as $L$ and all the black vertices as $R$.

\begin{tcolorbox}[colback=yellow!5!white,colframe=yellow!75!black,title=\textbf{Example 4.3} Bipartite Graph of a Ferrers Board]
  \begin{center}
      $\young(\bullet\circ,\circ\bullet,\bullet,\circ)$
  \end{center}
  If we were to represent it as a bipartite graph, we can have the black squares as our left vertices and white squares as our right vertices. And then we can connect to the two sets with edges based on their adjacency in the Ferrers board. For our Ferrers board, this gives us:
  \begin{center}
  \begin{tikzpicture}
   \graph[nodes={draw, circle,fill=myblue}, radius=.5cm,empty nodes, branch down=1 cm, grow right sep=4cm] {subgraph I_nm [V={a, b, c}, W={1,...,3}];
  a -- {1,2};
  b -- {1,2};
  c -- {2,3};
};
\end{tikzpicture}
\end{center}
\end{tcolorbox}
Now that we have a bipartite graph, we can work with this. If you can find a perfect bipartite matching on this bipartite graph then it will correlate with a certain tiling of the Ferrers board. This is because each edge in the bipartite graph corresponds to a potential domino being placed since each domino is just 2 adjacent squares. So if we can find a perfect bipartite matching, that correlates to disjoint dominoes that tile the entire board.\\\\
We showed that a bipartite graph has a perfect bipartite matching if and only if it follows Hall's theorem. So let's see if we can apply our Ferrers boards to Hall's theorem.\\\\
In order for Hall's theorem to hold, we need that all subgraphs of our board to follow the inequality $|N(S)| \geq |S|$.

\begin{tcolorbox}[enhanced,attach boxed title to top center={yshift=-3mm,yshifttext=-1mm},
colback=blue!5!white,colframe=blue!75!black,colbacktitle=red!80!black,
  title=Exercise 4.x,fonttitle=\bfseries,
  boxed title style={size=small,colframe=red!50!black} ]
\textbf{4.3 ($\star$)} Show that all bipartite graphs formed from Ferrers boards follow Hall's theorem.
\end{tcolorbox}
Once you show this fact, then it shows that all Ferrers boards with the white/blackness property is necessarily and sufficently tileable. Which answers our question.

\newpage
\section{\textbf{Generating Tilings}} This section deals with some more computer science specific parts of this problem that I looked into while solving this problem. I won't go into too much detail, but I will list things that I read in references if you are interested. Primarily it will talk about how you can use a computer to generate all tilings of a Ferrers board, or just tilings of objects in general. I won't go into too much detail for any of the algorithms, so if you're interested you can google and read for yourself.
\subsection{Max Flow}
One solution to find a spefic tiling is applying max flow algorithms. Max flow is a class of problems in combinatorial optimization which deals with finding a feasible flow through a network that gives you the maximum amount of flow.\\\\
One specific application of max flow is being able to find a perfect bipartite matching in a given bipartite graph. The most common algorithm to do so is the Ford-Fulkerson Algorithm.\\\\
I will quickly talk about the specific case for our Ferrers boards. We've already shown that we can create a bipartite graph out of our Ferrers boards that have an equal amount of black and white tiles. We take this graph and we add $2$ more nodes called the sink and the target. The sink will be connected to all of the left nodes while the target will be connected to all of the right nodes. Which will look something like this:
\begin{center}
  \includegraphics[width=100mm]{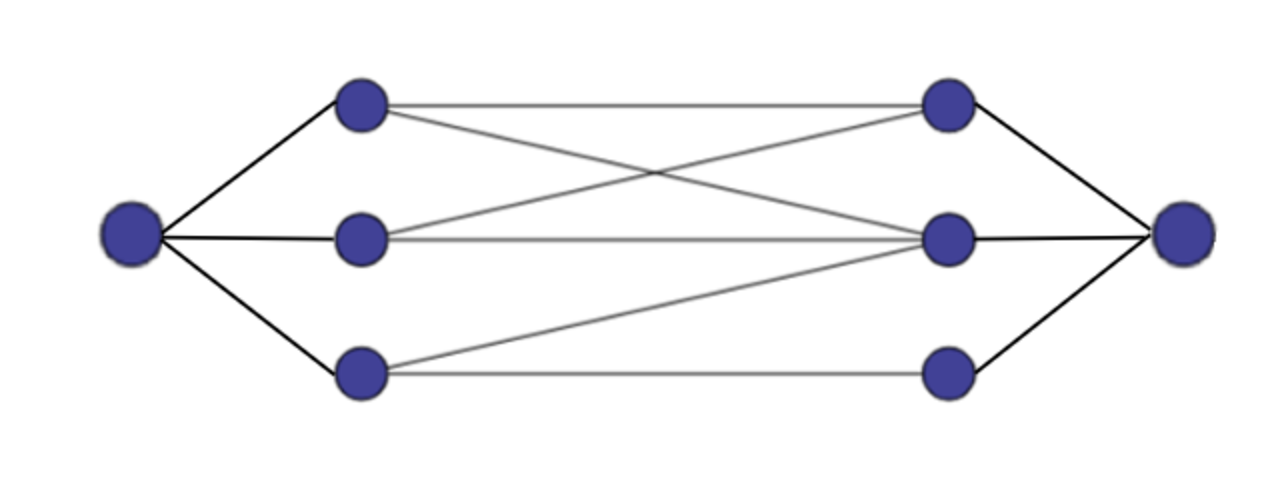}
\end{center}
We will also make this graph a directed graph, where all edges are pointing from the left to the right. So essentially all nodes from the sink will be flowing into the target. We will also make it so that each edge has a "capacity" of $1$. This will make more sense shortly.\\\\
We have now converted our tiling problem into a max flow problem, because finding the max flow on this graph now is the same as finding a perfect bipartite matching, which we showed is the same as finding a tiling. We can accomplish this with an algorithm called the Ford-Fulkerson Algorithm.

\begin{tcolorbox}
[colback=blue!5!white,colframe=blue!75!black,title=\textbf{Definition 5.1} (Ford-Fulkerson)]
The basic idea of the Ford-Fulkerson algorithm is to start with an empty matching M, and iteratively augment it until we can no longer increase the size of the matching. The algorithm works as follows:
\begin{enumerate}
    \item Start with an empty matching M.
    \item Find an augmenting path in the graph using a graph traversal algorithm, such as BFS or DFS. An augmenting path is a path from a vertex in U that is not matched to a vertex in V, to a vertex in V that is not matched, such that the path alternates between unmatched and matched edges.
    \item If an augmenting path is found, update the matching M by flipping the matching status of the edges along the augmenting path. This means that unmatched edges become matched, and matched edges become unmatched. This increases the size of the matching by 1.
    \item If no augmenting path is found, stop. The current matching M is a maximum matching.
\end{enumerate}
Once we have our M, that is equivalent to our perfect bipartite matching, and gives us our tiling.
\end{tcolorbox}
I would show an example, but I don't want to LaTeX it unfortunately so I'm going to leave it at that.
\subsection{Other Algorithms} 
\begin{enumerate}
    \item The Fischer-Kastelyn-Temperley Algorithm or FKT Algorithm, is an algorithm that allows us to count the number of perfect bipartite matchings in a graph in polynomial time. So although this won't allow us to generate our tilings, it will allow us to quickly count how many we have with the help of a computer. 
    \item The Fukuda-Matsui algorithm allows us to find all perfect bipartite matchings in $O(c(n+m)+n^{2.5})$ (where n denotes the number of vertices, m denotes the number of edges, and c denotes the number of perfect matchings in the given bipartite graph).
\end{enumerate}

\end{document}